\documentclass[12pt,reqno]{amsart}

\usepackage{amsmath}
\usepackage{amssymb}
\usepackage{amscd}
\usepackage{amsthm}
\usepackage{amsfonts}
\usepackage[utf8]{inputenc}
\usepackage[T1]{fontenc}
\usepackage{lmodern}
\usepackage{paralist}
\usepackage{epic}
\usepackage{enumitem}

\newtheorem{theorem}{Theorem}
\newtheorem{theorema}{Theorem}
\newtheorem{theoremb}{Theorem}
\newtheorem{theoremc}{Theorem}
\newtheorem{theoremd}{Theorem}
\newtheorem{theoreme}{Theorem}
\newtheorem{rem}[theorema]{Remark}
\newtheorem{lem}[theoremb]{Lemma}
\newtheorem{prop}[theoremc]{Proposition}
\newtheorem{cor}[theoremd]{Corollary}
\newtheorem{defi}[theoreme]{Definition}
\newtheorem{ex}[theoreme]{Example}

\newcommand\stab{\op{Stab}}
\newcommand\ann{\op{Ann}}

\newcommand\C{{\mathbb C}}

\newcommand\g{\mathfrak{g}}

\newcommand\h{\mathfrak{h}}
\renewcommand\b{\mathfrak{b}}
\renewcommand\H{\mathbb{H}}

\newcommand\op[1]{\mathop{\rm #1}\nolimits}

\newcommand\p{\mathfrak{p}}
\newcommand\R{{\mathbb R}}
\renewcommand\sl{\mathfrak{sl}}

\newcommand{\Vtm}{\mathbb{V}_{\theta_{\text{max}}}}

\begin{document}

\title[Minimal Projective Orbits of Semi-simple Lie Groups]{Minimal Projective Orbits of Semi-simple Lie Groups}
\author[H.\ Winther]{Henrik Winther}
\address{$\ddag$ Department of Mathematics and Statistics, Faculty of Science and Technology, UiT – The Arctic University of Norway, N-9037 Tromsø, Forskningsparken 1, B4, Norway.
E-mail: \textsc{henrik.winther@uit.no}.}
\maketitle

\begin{abstract}
	Let $G$ be a Lie group $G$ with representation $\rho$ on a real simple $G$-module $\mathbb{V}$. We will call the orbits of the induced action of $\rho$ on the projectivization $P\mathbb{V}$ the projective orbits, and projective orbits of lowest possible dimension will be called \emph{minimal}. 
	We show that when $G$ is semi-simple and non-compact, there exists a compact subgroup $K\subset G$ such that the minimal orbits of $G$ are in bijection with the minimal $K$-orbits on a $K$-invariant proper subspace $\mathbb{W}\subset \mathbb{V}$. 
	In the case that $G$ is split-real, $K$ is the trivial subgroup and there is a unique closed projective orbit, which is moreover of minimal dimension.
\end{abstract}

\section{Introduction}\label{Intro}
Let $G$ be a real semi-simple Lie group, with representation $\rho$ on a real non-trivial simple $G$-module $\mathbb{V}$. 
We are interested in the problem of classifying the minimal orbits of $G$ in $\mathbb{V}$, that is, the non-zero orbits of lowest possible dimension. 

The general solution is understood when the Lie group $G$ is complex and acting by a complex-linear representation. 
In that case, the minimal projective orbit is the orbit of the class of a highest weight vector, by the Borel theorem \cite{B}.
We will generalize that result in this note.

The related case of closed orbits in $\mathbb{V}$ can be understood by the method of minimal vectors. The following cases fall under this regime:  
\begin{itemize}
	\item complex algebraic groups acting by complex representations \cite{KN},
	\item real algebraic groups acting by real rational representations \cite{RS},
	\item real reductive Lie groups acting by linear representations \cite{BL}.
\end{itemize}
However, these methods do not give an immediate way to see the dimension of the orbits, and we emphasize that closed projective orbits need not lift to closed orbits in $\mathbb{V}$. 

The special case of the group $G=Sp(1)GL(n,\H)$, acting on the curvature- and torsion-modules of the corresponding $G$-structure, was considered in \cite{KWZ}. 
The current work is a generalization of the techniques developed in that paper.

The problem described above occurs frequently in differential geometry, for example when classifying geometric structures with large (local) symmetry algebras.
Suppose $M$ is a Cartan geometry modeled on $F/G$ which admits a non-trivial automorphism group $\text{Aut}(M)$.
The differential invariants of such a structure are typically sections of some associated vector bundle, $T \in \Gamma(E\rightarrow M)$. 
At a point $x\in M$, the image $H$ of the isotropy representation of the stabilizer subgroup $\text{Stab}(x)$ of that point is a linear Lie group. 
This linear group $H$ acts on the associated $G$-module $\mathbb{E}$ by restriction, and it must necessarily preserve the pointwise value $T_x\in E_x \simeq \mathbb{E}$ of all invariants. 
The linear group $H$ is then contained in the stabilizer subgroup of the orbit $G\cdot T_x$.

Thus, non-flat geometries with abundant automorphisms are related to orbits of small dimension in $\mathbb{E}$, but one might equivalently consider the orbits in the projectivization $P(\mathbb{E})$.

The dimension of an orbit $G\cdot v$ for  $v\in \mathbb{V}$ is the codimension of the stabilizer $H=\stab{v}$ in $G$. Therefore, minimizing the dimension of an orbit is the same as maximizing the dimension of the corresponding stabilizer.
Hence, an alternative approach to the one we will describe is to enumerate the subalgebras $\h$ of the Lie algebra $\g$ of $G$, and check whether each conjugacy class admits an invariant vector, starting with the highest dimensional subalgebras and checking successively smaller ones until an invariant vector is found.
The main computational benefit of our proposed solution to the problem of finding minimal orbits is that we avoid this process of enumerating subalgebras completely.

\subsection*{Acknowledgments} We would like to thank A. C. Ferreira and B. Kruglikov for useful discussion, and D. The for proposing an improved formulation of Corollary \ref{cor:realhighestweight}.
Further, we gratefully acknowledge the support of the Czech Science Foundation via the project (GA\v{C}R) No.\ GX19-28628X.
\section{Minimal orbits}
Let $G$ be a real, semi-simple, connected Lie group and $\mathbb{V}$ be a simple finite-dimensional $G$-module. 
Let $\g_{-k} \oplus \dots \oplus \g_{-1} \oplus \g_0 \oplus \g_1 \oplus \dots \oplus \g_k$ be a non-trivial $\mathbb{Z}$-grading on its Lie algebra $\g$, where $\g_i = 0$ for $|i|>k$, let $\p = \oplus_{i\ge 0}\g_i$ be the non-negatively graded (parabolic) subalgebra, and let $\p_+$ be the strictly positively graded nilpotent subalgebra.
For a non-trivial $\mathbb{Z}$-grading of a semi-simple Lie algebra, there exists a \emph{grading element} $Z\in \g_0$, which acts on each component $\g_i$ by the eigenvalue $i$ in the adjoint representation \cite{capslovak}. 
In fact, $Z$ is a non-compact element of a maximally non-compact Cartan subalgebra, and so it is diagonalizable over $\mathbb{R}$ for any $\g$-representation. 
Let $\bigoplus_{\theta\in \{\theta_{\text{min}}, \dots , \theta_{\text{max}}\}} \mathbb{V}_{\theta}$ be the eigenspace decomposition of $\mathbb{V}$ with respect to $Z$.
From now on we shall refer to the $G$-orbits in $\mathbb{V}$ as \emph{linear orbits} and the $G$-orbits in $P\mathbb{V}$ as \emph{projective orbits}.

\begin{prop}\label{prop:gradedaction}
	Let $g_i \in \g_i$. Then, for $v\in \mathbb{V}_{\theta}$, we have $g_i(v) \in \mathbb{V}_{\theta+i}$.
\end{prop}
\begin{proof}
	By definition, $[Z,g_i] = i g_i$. Thus, applying each side of the equality to the vector $v$, we get $Z\cdot g_i(v) = (\theta + i)g_i(v)$.
\end{proof}
\begin{defi}
	Let $\mathbb{V}$ be a $\g$-module and $\rho$ be the associated $\g$--representation. Then for a subalgebra $\h \in \g$, we define the kernel $\ker \h = \cap_{h\in\h}\ker(\rho(h))$.
\end{defi}

\begin{prop}\label{prop:g0andkerp}
	We have $\ker \p_+ = \ker \g_1 = \mathbb{V}_{\theta_{max}}$, and moreover the connected linear Lie group $G_0$ generated by exponentials of $\g_0$ acts irreducibly by its restricted representation to $\mathbb{V}_{\theta_{max}}$.
\end{prop}
\begin{proof}
	Let $v,w \in \mathbb{V}$ be non-zero. Then, since $\mathbb{V}$ is simple with respect to $G$, and $G$ is assumed connected, the vector $w$ generates $\mathbb{V}$ over $\g$, and we may write 
	\begin{equation}
		v = \sum_{i=1}^k \prod_{j=1}^{l_i} g^{(i,j)}_{\alpha_{ij}} w,
		\label{eq:simpleactiongenerate}
	\end{equation}
	where $g^{(i,j)}_{\alpha_{ij}}\in \g_{\alpha_{ij}}$, and the products are ordered. The products can be understood as compositions of linear operators representing the module action of $g^{(i,j)}_{\alpha_{i,j}}$, or in the universal enveloping algebra $\mathfrak{U}(\g)$.
	In fact, we may assume that the products are ordered by (non-strictly) increasing order of gradation (left to right), as for $g_i \in \g_i, g_j \in \g_j $ with $i>j$, we can rewrite $g_{i}g_j = g_j g_i + g_{i+j}^\prime$ for $g_{i+j}^\prime = [g_i, g_j]$, and split the term in two.
	Now if $w\in \ker \p_+$, then each non-zero term in (\ref{eq:simpleactiongenerate}) has only factors with non-positive gradation, as all positively graded factors annihilate $w$. 
	Hence, setting $v \in \mathbb{V}_{\theta_{\text{max}}}$, we get that each non-zero term must contain only factors of gradation $0$, by Prop. \ref{prop:gradedaction}.
	Therefore, $\ker \p_+ = \ker\g_1 = \mathbb{V}_{\theta_{\text{max}}}$, and $\g_0$ acts irreducibly by its restricted representation to $\mathbb{V}_{\theta_{\text{max}}}$, hence the same is true for $G_0$.
\end{proof}
The following result is a small modification of \cite[Lemma 3]{KWZ}.
\begin{lem}\label{lem:existenceminimal}
Let $\mathbb{V}$ be a real, finite-dimensional, simple $G$-module for a real connected Lie group $G$, such that the center of $G$ acts on $\mathbb{V}$ either by real scalars, or by complex scalars. In the latter case, $\mathbb{V}$ is realized as the underlying real module of a complex module. Then there exists a projective orbit $G\cdot [v] \subset P\mathbb{V}$ which is closed and of minimal dimension. Moreover, any projective orbit of minimal dimension is closed.
\end{lem}

\begin{proof}
	We may replace the group by its image under the associated representation, since this will not affect the (projective) orbits. 
	Therefore, we may assume without loss of generality that the representation is faithful.
	Because the representation is irreducible, $\mathfrak{g}$ is the direct sum of a semi-simple ideal and a central subalgebra \cite[Chap. 3]{J}. 
	Our assumption on the center action means that the module $\mathbb{V}$ is a real or complex tensor product of a (real or complex)-irreducible module $\mathbb{V}_0$ over the semi-simple ideal and a (real or complex) one-dimensional module $\R$ or $\C$ over the center. 
	If the tensor product is complex, then it is taken with respect to some invariant complex structures on the underlying real modules.
	Let $G \cdot [v]\subset P\mathbb{V}$ be a projective $G$--orbit of minimal dimension $d$.
	We want to prove it is closed.
	Consider the complexifications of the group, the action and the module.
	The element $v \in \mathbb{V}+0\cdot i \subset \mathbb{V}^\C$ determines the complex orbit $G^\mathbb{C} \cdot [v] \subset P \mathbb{V}^\C$ of \emph{complex} dimension $d$, since the annihilator of $v$ in $\g^\C$ intersects $\g$ by the annihilator of $v$ in $\g$.
	If the closure of the orbit $G \cdot [v]$ contains another orbit $G\cdot [v^\prime]$ (necessarily of the same dimension $d$), then the closure of the complex orbit $G^\C \cdot [v]$ contains the complex orbit $G^\C \cdot [v^\prime]$ (again of the same complex dimension $d$).
	To exclude the latter, note that the action of $G^\C$ on $P\mathbb{V}^\C$ is algebraic, because both the semi-simple action on $\mathbb{V}_0^\C$ and the central action on $\C$ are algebraic (as we replaced the group by its image) \cite[Chap. 3]{OV}. 
	But then the boundary of any orbit can only contain orbits of strictly lower dimensions, which are less than $d$ (see the proof of \cite[III, \S 1.5]{OV}).
	This contradiction proves the claim.
\end{proof}

\begin{theorem}\label{thm:closedorbitkerp}
	Let $G$ be a real, semi-simple, connected Lie group and $\mathbb{V}$ be a simple finite-dimensional $G$-module. 
	Let $\g_{-k} \oplus \dots \oplus \g_{-1} \oplus \g_0 \oplus \g_1 \oplus \dots \oplus \g_k$ be a non-trivial $\mathbb{Z}$-grading on its Lie algebra $\g$, with $\p = \oplus_{i\ge 0}\g_i$ its non-negatively graded (parabolic) subalgebra, and $\p_+$ the strictly positively graded nilpotent subalgebra.
	Suppose that the orbit $G\cdot [v]\subset P \mathbb{V}$ is closed. 
	Then there exists a nonzero element $w \in \ker \p_+$ such that $[w]\in G \cdot [v]$.
\end{theorem}

\begin{proof}
	Any linear orbit contains some element $v=\sum_\theta v_{\theta}$ with nonzero component $v_{\theta_{\text{max}}}$ in $\mathbb{V}_{\theta_\text{max}}$ (or any other $Z$-eigenspace), since otherwise the orbit would span a proper submodule. 
	So we may assume without loss of generality that $v_{\theta_{\text{max}}}\not= 0$.
	Consider the sequence $w_m, m\in \mathbb{N}$ given by $w_m = \exp(tZ)|_{t=m}(v)$. We have $w_m = \sum_\theta e^{m\theta}v_{\theta}$ in the $Z$-eigen-decomposition. The sequence $[w_m]$ converges to $[v_{\theta_{\max}}]$ in $P\mathbb{V}$, and since $G\cdot [v]$ is assumed to be closed, the limit is contained in the same projective orbit. Hence $[v_{\theta_{\text{max}}}]\in G \cdot [v]$, so we let $w = v_{\theta_{\text{max}}}$, and obtain the result.
\end{proof}

\begin{cor}\label{cor:realhighestweight}
	Let $G$ be split-real and $\mathbb{V}$ a simple finite-dimensional $G$-module. Then there is a unique closed projective orbit, and it is generated by a real highest weight vector. This orbit is of minimal dimension.
\end{cor}
\begin{proof}
	Since $G$ is split, there exists a totally non-compact Cartan subalgebra $\h\subset\g$, and this admits a $\mathbb{Z}$-grading corresponding to a real slice of the Borel subalgebra in the complexification $\g^\C$. 
	In particular, we have $\h = \g_0$, abelian and real diagonalizable.
	By Prop. \ref{prop:g0andkerp}, $\ker \p_+ = \mathbb{V}_{\theta_{\text{max}}}$ is simple, so it is one-dimensional. 
	Thus $\mathbb{V}$ is a real highest weight representation generated by $v\in\mathbb{V}_{\theta_{\text{max}}}$, and by Lemma \ref{lem:existenceminimal}, a closed minimal orbit exists. Therefore, by Theorem \ref{thm:closedorbitkerp}, the unique closed projective orbit is $G\cdot [v]$.
\end{proof}
Let us now consider the $\mathbb{Z}$--grading associated to the \emph{minimal parabolic subalgebra} $\b\subset \g$. 
When $\g$ is split-real, this is a real form of the Borel subalgebra, but in general this does not hold. 
Instead it is the parabolic subalgebra associated to the Satake diagram of $\g$ with all white nodes crossed, therefore we have that ${\b_0}_{ss} = \mathfrak{k}_1 \oplus \dots \oplus \mathfrak{k}_l$ is the direct sum of a finite number of compact subalgebras. 
\begin{prop}\label{prop:thegoodparabolic}
	Let $G$ be a non-compact semi-simple Lie group and $\mathbb{V}$ a simple, finite-dimensional and non-trivial $G$--module. There exists a non-trivial $\mathbb{Z}$-grading on $\g$ such that for any any closed projective orbit, there exists a representative $w\in\mathbb{V}_{\theta_{\text{max}}}$ such that $\ann(w)\subset\g_0\oplus\g_+$.
\end{prop}
\begin{proof}
	Let $\b$ be the minimal parabolic subalgebra, and $\sigma$ a Cartan involution of $\g$ such that the Cartan subalgebra contained in $\b_0$ is $\sigma$-stable.
	Let us first fix a minimal orbit with representative $w \in \Vtm$, and $\b_+\subset \ann(w)$, via Theorem \ref{thm:closedorbitkerp}.
	Suppose we have $x_- \in \ann(w) \cap \g_-$ nonzero.
	Then we get $x_+ = \sigma x_- \in \b_+$, and $x_0 = [x_-, x_+] \in \b_0$.
	Moreover, since $\sigma x_0 = [\sigma x_-,\sigma x_+] = [x_+,x_-] = -x_0$, we have that $x_0$ is a non-compact Cartan element.
	Now $\langle x_-, x_0, x_+ \rangle \simeq \sl_2(\R)$, with the elements giving a standard basis.
	The subspace $\Vtm$ is invariant with respect to this $\sl_2$, in particular it is trivial, since it consists of highest weight vectors (as $x_+ \in \b_+$) with weight zero with respect to its Cartan subalgebra spanned by $x_0$.
	The latter follows because noncompact Cartan elements of $\g$ act by real scalars due to Prop. \ref{prop:g0andkerp}, and $x_0$ has the eigenvector $w$ with eigenvalue $0$.
	Thus if an element of $\g_-$ annihilates $w$, then it also annihilates the representatives of all other closed projective orbits.
	Let $\p = \langle \b, \ann(w) \rangle$.
	This is a proper subalgebra because $\mathbb{V}$ is a non-trivial simple module, and it is parabolic because it contains $\b$. The  $\mathbb{Z}$--grading associated to $\p$ then satisfies $\ker(\b_+) \subset \ker(\p_+)$, so for each closed projective orbit, there exists a representative $w^\prime\in \Vtm \subset \ker(\p_+)$ for which we have $\ann(w^\prime)\subset\p_0\oplus\p_+$.
\end{proof}

\begin{theorem}\label{thm:Korbits}
	Let $G$ be a real semi-simple and non-compact Lie group, and $\mathbb{V}$ a real finite-dimensional simple $G$--module.
	Then there exists a compact subgroup $K\subset G$ such that the minimal projective $G$-orbits are in bijective correspondence with minimal projective $K$-orbits of a simple $K$-submodule $\mathbb{W} \subset \mathbb{V}$.
\end{theorem}
\begin{proof}
	We define the following operation on pairs $(G,\mathbb{V})$, where $G$ is a real connected semi-simple Lie group and $\mathbb{V}$ is a simple $G$-module: $(G,\mathbb{V}) \mapsto (G^\prime, \mathbb{V}^\prime)$, where we apply Prop \ref{prop:thegoodparabolic} to $(G,\mathbb{V})$, to obtain a $\mathbb{Z}$--grading on $\g$, and let $G^\prime$ be the semi-simple Lie group generated by those simple ideals of $\g_0$ which act non-trivially on $\mathbb{V}_{\theta_{\text{max}}}$, and $\mathbb{V}^\prime = \mathbb{V}_{\theta_{\text{max}}}$.
	Next we generate a sequence of pairs $(G^m,\mathbb{V}^m)$, by setting $G^{m+1} = (G^m)^\prime$, $\mathbb{V}^{m+1} = (\mathbb{V}^m)^\prime$, starting with $G^0 = G, \mathbb{V}^0 = \mathbb{V}$.
	Note that by Prop. \ref{prop:g0andkerp}, each $\mathbb{V}^m$ is simple, satisfying one of the assumptions of Prop. \ref{prop:thegoodparabolic}. 
	However the sequence cannot be continued infinitely, as the procedure must end at some step $k$, which is the first step for which $G^k$ is a compact group (which might be the trivial subgroup \{e\}), for which Prop. \ref{prop:thegoodparabolic} can not be applied.
	For any $\mathbb{Z}$--grading, the center of $G_0$ acts via real (or complex) scalars on $\mathbb{V}_{\theta_{\text{max}}}$, and a given central element is either in the kernel of the restricted representation, or it admits no invariant vector.
	Therefore, the center may be disregarded for the purposes of finding minimal orbits at each step.
	Then a combination of Theorem \ref{thm:closedorbitkerp} and the properties guaranteed by Prop. \ref{prop:thegoodparabolic} gives that the minimal projective orbits are in bijection with each other at each step (in fact, the representatives are the same elements of $\mathbb{V}$).
	Thus, set $K = G^k$, $\mathbb{W} = \mathbb{V}^k$.
	This yields the result.
\end{proof}
\begin{rem}
	By using the proof of Prop. \ref{prop:thegoodparabolic}, it is possible to see the Lie algebra of the subgroup $K$ from Theorem \ref{thm:Korbits} via the Satake diagram of $\g$, decorated with the coefficients of the highest weight of the complex simple module associated to $\mathbb{V}$, in the fundamental weight basis (as in \cite{capslovak}). 
	The following rule can be used: A white node should be crossed if either it has a non-zero coefficient, or it is immediately adjacent to a fully black connected sub-diagram which has at least one non-zero coefficient.
	Removing the crossed nodes yields some number of Satake diagrams decorated with all zeroes, which may be disregarded, and some completely black Satake diagrams, each decorated with at least one non-zero coefficient. The latter form the Lie algebra of $K$, and the remaining decorations give the highest weight of the complex simple module associated to $\mathbb{W}$, in the fundamental weight basis.
	\label{rem:easymode}
\end{rem}

\begin{ex}\label{ex:quatcurvature}
	The results \cite[Corollary 4.6, Corollary 4.8]{KWZ} can be seen retrospectively as an application of Theorem \ref{thm:Korbits}. 
	There, a minimal projective orbit of the Lie group $Sp(1)SL(n,\H)$ is computed for both the curvature module $\mathbb{V}^{II} = S^3_\C \H^{n\ast} \odot \H^n$ and the torsion module $\mathbb{V}^{I} =S^3_\C \H \otimes \Lambda^2_\C \H^{n\ast} \odot \H^n$.
	In particular, the Satake diagrams we described in Remark \ref{rem:easymode} can be seen in \cite[Remark 5]{KWZ}, and we reproduce them here:
	\begin{align*}
	& \\
	&
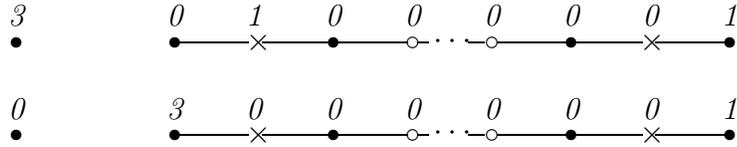
\begin{picture}(0,0)
	\put(-160,0){\circle*{4}}
	\put(-100,0){\circle*{4}}
	\put(-100,0){\line(10,0){28}}
	\put(-73,-3){$\times$}
	\put(-67,0){\line(10,0){25}}
	\put(-40,0){\circle*{4}}
	\put(-40,0){\line(10,0){28}}
	\put(-10,0){\circle{4}}
	\put(-8,0){\line(10,0){4}}
	\put(-2,-2){$\cdots$}
	\put(11,0){\line(10,0){6}}
	\put(20,0){\circle{4}}
	\put(22,0){\line(10,0){28}}
	\put(50,0){\circle*{4}}
	\put(50,0){\line(10,0){27}}
	\put(76,-3){$\times$}
	\put(82,0){\line(10,0){28}}
	\put(110,0){\circle*{4}}
	\multiputlist(-160,10)(30,0){3,,0,1,0,0}
	\multiputlist(20,10)(30,0){0,0,0,1}
	\end{picture} \\
	& \\
	&\begin{picture}(0,0)
	\put(-160,0){\circle*{4}}
	\put(-100,0){\circle*{4}}
	\put(-100,0){\line(10,0){28}}
	\put(-73,-3){$\times$}
	\put(-67,0){\line(10,0){25}}
	\put(-40,0){\circle*{4}}
	\put(-40,0){\line(10,0){28}}
	\put(-10,0){\circle{4}}
	\put(-8,0){\line(10,0){4}}
	\put(-2,-2){$\cdots$}
	\put(11,0){\line(10,0){6}}
	\put(20,0){\circle{4}}
	\put(22,0){\line(10,0){28}}
	\put(50,0){\circle*{4}}
	\put(50,0){\line(10,0){27}}
	\put(76,-3){$\times$}
	\put(82,0){\line(10,0){28}}
	\put(110,0){\circle*{4}}
	\multiputlist(-160,10)(30,0){0,,3,0,0,0}
	\multiputlist(20,10)(30,0){0,0,0,1}
	\end{picture}\\
	&
	\end{align*}
	Thus, on removing the crossed nodes, we obtain
	\begin{align*}
	& \\
	&
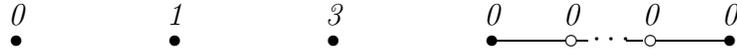
\begin{picture}(0,0)
	\put(-160,0){\circle*{4}}
	\put(-100,0){\circle*{4}}
	\put(-40,0){\circle*{4}}
	\put(20,0){\circle*{4}}
	\put(20,0){\line(10,0){28}}
	\put(50,0){\circle{4}}
	\put(52,0){\line(10,0){4}}
	\put(58,-2){$\cdots$}
	\put(71,0){\line(10,0){6}}
	\put(80,0){\circle{4}}
	\put(82,0){\line(10,0){28}}
	\put(110,0){\circle*{4}}
	\multiputlist(-160,10)(30,0){0,,1,,3,,0,0,0,0}
	\end{picture}\\
	&
	\end{align*}
	This shows that in both cases, the minimal projective orbits are in bijective correspondence with the minimal projective $K=Sp(1)Sp(1)$--orbits on the simple $K$--module $S^3_\C\H \otimes_\C \H$, which is of real dimension 8.
\end{ex}
\begin{rem}\label{rem:dimstability}
	It is worth noting the following phenomenon illustrated by Example \ref{ex:quatcurvature}: for parametrized families of groups (e.g. $Sp(1)SL(n,\H)$ with parameter $n$), with representations expressed in a fixed way in terms of weights, the $K$--module $\mathbb{W}$, and hence the set of minimal projective orbits, does not necessarily depend on the parameter. 
\end{rem}

\begin{ex}
	Consider $SO(p,q)$ for $0\le p<q$, and the module $\mathbb{V} = \Lambda^{p+1}(\R^{p+q})$.
	This admits a unique closed (minimal) projective orbit, because $K = SO(q-p)$ and $\mathbb{W}$ is the standard representation of $K$, which is sphere-transitive.
        A representative can be constructed by writing $\R^{p+q} = \R^{p,p} \oplus \R^{q-p}$. 
	Then taking a null plane $\Pi^p\subset \R^{p,p}$, and $v\in \R^{q-p}$, we get the unique closed projective orbit $G\cdot [\Lambda^p\Pi^p \wedge v]$. 
\end{ex}

\begin{rem}
	More examples with unique closed projective orbit can be constructed by taking $\mathbb{W}$ to be a sphere-transitive representation of a simple compact Lie  group $K$, and extending its Satake diagram to that of a non-compact simple Lie algebra.
	The eligible groups $K$ are: $SO(n)$, $SU(n)$, $Sp(n)Sp(1)$, $Sp(n)$, $Spin(7)$ and $Spin(9)$, where in the latter two cases $\mathbb{W}$ is the spin module (see \cite{MS}).
	The Satake diagram of $G_2$ does not extend to a Satake diagram of a noncompact Lie algebra.
\end{rem}


\end{document}